\documentstyle[12pt,amssymb]{amsart}
\newcommand{\SH}{{\goth H}}
\newcommand{\Th}{\Theta}
\newcommand{\lan}{\langle}
\newcommand{\ran}{\rangle}

\newcommand{\SL}{\operatorname{SL}}
\renewcommand{\Sp}{\operatorname{Sp}}

\newcommand{\th}{\theta}

\newcommand{\MM}{{\cal M}}

\newcommand{\CC}{{\cal C}}
\newcommand{\Ga}{\Gamma}
\newcommand{\ga}{\gamma}
\newcommand{\si}{\sigma}

\newcommand{\rk}{\operatorname{rk}}

\newcommand{\ed}{\qed\vspace{3mm}}

\newcommand{\Spec}{\operatorname{Spec}}

\newtheorem{thm}{Theorem}[section]
\newtheorem{prop}[thm]{Proposition}
\newtheorem{lem}[thm]{Lemma}
\newtheorem{cor}[thm]{Corollary}

\theoremstyle{definition}

\newenvironment{rem}{\vspace{3mm}
\noindent {\bf Remark.}}{\vspace{3mm}}
\newenvironment{rems}{\vspace{3mm}
\noindent {\bf Remarks.}}{\vspace{3mm}}

\numberwithin{equation}{section}
                
\newcommand{\Pf}{\noindent {\it Proof}}
\newcommand{\Pic}{\operatorname{Pic}}

\renewcommand{\a}{\alpha}
\renewcommand{\b}{\beta}
\renewcommand{\S}{{\bold S}}

\newcommand{\wt}{\widetilde}
\renewcommand{\mod}{\operatorname{mod}}

\renewcommand{\AA}{{\cal A}}

\newcommand{\FF}{{\cal F}}

\newcommand{\Z}{{\Bbb Z}}

\newcommand{\la}{\lambda}
\newcommand{\PP}{\cal P}

\newcommand{\PGL}{\operatorname{PGL}}

\newcommand{\ov}{\overline}

\newcommand{\om}{\omega}
\newcommand{\ra}{\rightarrow}

\newcommand{\id}{\operatorname{id}}
\newcommand{\G}{{\Bbb G}}

\newcommand{\ot}{\otimes}

\newcommand{\We}{\bigwedge}
\renewcommand{\O}{{\cal O}}
\renewcommand{\P}{{\Bbb P}}

\newcommand{\du}{\vee}

\newcommand{\de}{\delta}

\newcommand{\D}{{\cal D}}
\newcommand{\De}{\Delta}

\newcommand{\sub}{\subset}
\newcommand{\Nm}{\operatorname{N}}

\title{Determinant bundles for abelian schemes}
\author{A. Polishchuk}
                
\begin{document}
        
        
\maketitle

\begin{abstract} To a symmetric, relatively ample line bundle
on an abelian scheme one can associate  a linear
combination of the determinant bundle and the
relative canonical bundle, which is a torsion
element in the Picard group of the base. We improve the
bound on the order of this element found by Faltings and Chai.
In particular, we obtain an optimal bound when
the degree of the line bundle $d$ is odd and
the set  of  residue  characteristics  of  the  base  does   not
intersect the set of primes $p$ dividing $d$,
such that $p\equiv -1\mod(4)$ and $p\le 2g-1$,
where $g$ is the relative dimension of the abelian scheme.
Also, we show that in some cases these torsion elements generate
the entire torsion subgroup in the Picard group of the corresponding
moduli stack.
\end{abstract}

Let $L$ be a relatively ample line
bundle on an abelian scheme $\pi:A\ra S$,
trivialized along the zero section. Assume that
$L$ is symmetric, i.~e. $[-1]_A^*L\simeq L$.
We denote
by $\phi_L:A\ra\hat{A}$ the corresponding self-dual
homomorphism (where $\hat{A}$ is the dual abelian scheme).
Let $d=\rk\pi_*L$, so that $d^2$ is the degree of $\phi_L$.
Then Faltings and Chai proved in \cite{FC}, I, 5.1
the following equality in $\Pic(S)$:
$$8\cdot d^3\cdot\det(\pi_*L)=-4\cdot d^4\cdot\ov{\om}_A$$
where we denote by $\ov{\om}_A$ the restriction of the
relative canonical bundle $\om_{A/S}$ to the zero section.
In other words, the element
$$\De(L):=2\cdot\det(\pi_*L)+d\cdot\ov{\om}_A$$
of $\Pic(S)$ is annihilated by $4d^3$. It is known from the
transformation theory of theta-functions (see \cite{MB2})
that this result is sharp for principal polarizations ($d=1$).
In the case of analitic families of complex abelian varieties
A.~Kouvidakis showed in
\cite{Kouv} using theta functions that
if the type of polarization is $(d_1,\ldots,d_g)$ with
$d_1|\ldots|d_g$ then  $4\cdot\De(L)=0$ except when
$3|d_g$ and $d_{g-1}\not\equiv 0\mod(3)$. In the latter
case he proved that $12\cdot\De(L)=0$.
This suggests that one can try to eliminate the factor $d^3$
in general situation. We prove that one can do this outside
certain set of prime divisors of $d$. In particular we
explain the appearence of factor $3$ above algebraically
(see Theorem \ref{isogmain}).
Here are the precise statements.
               
\begin{thm}\label{main1}
Let $L$ be a symmetric, relatively ample line
bundle over an abelian scheme $A/S$
of relative dimension $g$, trivialized  along  the
zero section. Then
\begin{enumerate}
\item $2^{n_2}\cdot d'\cdot\De(L)=0$ \hfill\\
where $d'=\prod p^{n_p}$ is the product
of powers of primes $p$ dividing $d$
such  that  $p\equiv  -1\mod(4)$ and $p\le \max(2g-1,3)$;
$n_p=1$ if $p\neq 3$ and $\frac{2g+1}{3}< p<2g-1$,
$n_p=2$ if $p=2g-1\neq 3$, otherwise  $n_p=v_p(d)$;
$n_2=2+3v_2(d)$.
\item There
exists an integer $N(g)>0$ depending only on $g$, such that
$$N(g)\cdot\De(L)=0.$$
\end{enumerate}
\end{thm}
                          
We  can  get  sharper  bounds  under some restrictions on the
residue characteristics of $S$. For every abelian group $K$
and a prime number $p$
we denote by $K^{(p)}$ the subgroup of elements of $K$
annihilated by some power of $p$.
Note that when $p$ is not among the residue characteristics of
$S$ we can define the $p$-type of polarization as the
type $(p^{n_1},\ldots,p^{n_g})$ of the finite symplectic group
$K(L)^{(p)}$, where $K(L)=\ker(\phi_L:A\ra\hat{A})$.
Here $n_1\le\ldots\le n_g$ are locally constant
functions on $S$.
              
\begin{thm}\label{main2}
Let $p$ be an odd prime     number,
$\De(L)^{(p)}\in\Pic(S)^{(p)}$  be the $p$-primary component
of $\De(L)$.
Let  $S[\frac{1}{p}]\sub  S$ be the open subscheme where $p$
is invertible. Then
$$\De(L)^{(p)}|_{S[\frac{1}{p}]}=0$$
in $\Pic(S[\frac{1}{p}])$ unless
$p=3$  and   the   $3$-type   of   the   polarization   over
$S[\frac{1}{3}]$ is $(1,\ldots,1,3^k)$, $k>0$. In the
latter case one has
$$3\cdot\De(L)^{(3)}|_{S[\frac{1}{3}]}=0$$
in $\Pic(S[\frac{1}{3}])$.
\end{thm}
             
\begin{rems} 1. In fact, the
equalities of Theorems \ref{main1} and \ref{main2}
can be realized by {\it canonical}
(i.e. compatible with arbitrary base changes)
isomorphisms of line bundles.
More precisely,
let $X_{g,d}$ be the moduli stack over $\Spec(\Z)$ classifying data
$(A,L)$ as above.
Since the construction of $\De(L)$ commutes
with arbitrary base changes, we can consider $\De(L)$ as an
element of the Picard group $\Pic(X_{g,d})$. Now we
claim that the equality of Theorem \ref{main1} holds in
$\Pic(X_{g,d})$ while the equality of Theorem \ref{main2}
holds in $\Pic(X_{g,d}[\frac{1}{p}])$.
This follows from the fact that there exists a
$\PGL_N$-torsor $\wt{X}_{g,d}$ over $X_{g,d}$ which is
represented by a scheme.
Indeed, $\wt{X}_{g,d}$ is obtained by adding a basis of
$\pi_*(L^3)$ (considered up to constant) to the  above  data
$(A,L)$. Then the representing scheme can be constructed as in
\cite{GIT} using Hilbert schemes. Now it suffices to prove the
triviality of the induced $\PGL_N$-equivariant line bundle on
$\wt{X}_{g,d}$. Applying  Theorem  \ref{main1}  (resp.  Theorem
\ref{main2}) to $\wt{X}_{g,d}$ we get some trivialization of
this line bundle. However, as $\PGL_N$ has no non-trivial
characters this trivialization is automatically compatible with
$\PGL_N$-action.
  
\noindent
2.  These  results  can be extended to line
bundles $L$ which are not necessarily ample,
but are {\it non-degenerate}, i.~e. such that
the corresponding homomorphism $\phi_L:A\ra\hat{A}$
is an isogeny. In this case there is a locally constant
function $i(L)$ on $S$ such that $R^i\pi_*(L)=0$ for
$i\neq i(L)$ and $R^{i(L)}\pi_*L$ is locally-free of rank
$d$. One has by definition
$\det\pi_*L:=(-1)^{i(L)}\det(R^{i(L)}\pi_*L)$. Also we should
replace $d$ by $(-1)^{i(L)}d$ when defining $\De(L)$:
$$\De(L):=2\cdot\det\pi_*L+(-1)^{i(L)}\cdot d\cdot\ov{\om}_A.$$
Then our  argument  goes  through  for  such  $L$.  Note
especially that in section \ref{meth} we rely heavily on the fact
that the function $n\mapsto\det\pi_*L^n$ has finite degree.
However, this is true for any line bundle ---
see Lemma \ref{degree}.
         
\noindent
3. When $d$ is even one can consider an element
$\De'(L)=\det\pi_*L+\frac{d}{2}\cdot\ov{\om}_A$
in $\Pic(S)$ such that $\De(L)=2\cdot\De'(L)$.
Kouvidakis proved in \cite{Kouv} that
for a totally symmetric line bundle $L$
one always has $3\cdot\De'(L)=0$ if $g\ge 3$ and the characteristic
is zero (furthermore, $\De'(L)=0$ if the $3$-polarization
type is not $(1,\ldots,1,3^k)$, $k>0$). It would be
nice to extend this result to the case of
positive   characteristics. So far, we only can control the behavior
of these elements under isogeny of odd degree
(see remark   after   Theorem \ref{isogmain2}).
Also, in the principally polarized case (i.~e. when $d=1$) one still
has $\De'(L^{2n})=0$ for $g\ge 3$ and
any scheme $S$ whose residue characteristics are prime to $2$.
This follows from the fact that the fibers of the relevant stacks
over $\Spec(\Z[\frac{1}{2}])$ are smooth and irreducible
(compare with the proof of Theorem \ref{torsion}).
\end{rems}
             
The following corollary describes the cases  when we get an
optimal result.
              
\begin{cor}  Assume that $d$ is odd and that every prime $p$ dividing $d$,
such  that  $p\equiv  -1\mod(4)$ and $p\le \max(2g-1,3)$,  is not
among residue characteristics of $S$.
Then one has $$4\cdot\De(L)=0$$
unless    the    $3$-type    of    the    polarization    is
$(1,\ldots,1,3^k)$, $k>0$. In the latter case one has
$$12\cdot\De(L)=0.$$
\end{cor}
        
Besides the idea of Faltings and Chai,
the crucial step in the proof of these theorems is
the relation between determinant bundles of $L$ and
$\a^*L$ where $\a$ is an isogeny of abelian schemes,
worked out in section \ref{behisog}.
Essentially this boils down to computation of norms
of symmetric line bundles with cube structures over
finite flat subgroups in abelian schemes. Theorem \ref{main1} is
proved in sections \ref{meth} and \ref{arsec}, while
Theorem \ref{main2} is proved in section \ref{charsec}.
In section \ref{comp} we calculate elements $\De(L)$ in the case $g=1$
and describe some linear relations between $\De(L^n)$ for different $n$
in higher dimensions. In particular, we prove that for $d=1$
and  $g\ge  2$ one has the following relation
$$\De(L^n)=n^{g-1}\cdot\De(L)$$
for any odd $n$.
Also  following  Mumford's approach we determine the torsion
subgroup in the Picard group of
the moduli stack $\wt{\AA}_g^+[\frac{1}{2}]$ of principally polarized
abelian varieties with even symmetric theta divisor (localized outside
characterstic $2$) for $g\ge 3$.
Recall that the divisor (and the corresponding line
bundle) is called even or odd depending on the parity of
its multiplicity at zero.
It turns out that the torsion subgroup in $\Pic(\wt{\AA}_g^+[\frac{1}{2}])$
is cyclic of order 4, hence it is generated by our element $\De(L)$.
           
\vspace{3mm}
               
\noindent
{\bf Notation}. Throughout this paper $S$ denotes the base
scheme and for any scheme over $S$ we denote by
$\pi$ its projection to $S$ (hoping that this will not lead
to confusion). For any morphism of schemes $f$ we
denote by $f_*$, $f^*$ and $f^!$
the corresponding standard functors between derived
categories of quasi-coherent sheaves (e.~g. $f_*$ denotes
the right derived functor of the push-forward functor).
For a vector bundle $V$ of rank $r$
we denote $\det V=\bigwedge^r(V)$. This definition can be
naturally extended to complexes (see \cite{De}).
For all abelian schemes over
$S$ we denote by $e$ the zero section. For an abelian scheme
$A$ (resp. morphism of abelian schemes $f$) we denote
by $\hat{A}$ (resp. $\hat{f}$) the dual abelian scheme
(resp. dual morphism). Also for every integer $n$ we denote by $[n]_A:A\ra A$
the multiplication by $n$ on $A$, and  by  $A_n\sub  A$  its
kernel.  For  every  line  bundle  $L$  on  $A$ we denote by
$\phi_L:A\ra\hat{A}$ the corresponding morphism of abelian
schemes (see \cite{MuAb}).
For a finite flat group scheme $H/S$ we denote by $|H|$
its order which is a locally constant function on $S$, so
for an integer $k$ the condition $|H|>k$ means that the order
of $H$ is greater than $k$ over every connected component of
$S$. We mostly use additive
notation for the group law in the Picard group.

{\it Acknowledgment}. I am grateful to E.~Goren, B.~Gross,
B.~Mazur, T.~Pantev and G.~Pappas for many helpful discussions.
       
\section{The behavior under isogenies}
\label{behisog}
                
In this section we study the relation between $\De(L)$
and $\De(\a^*L)$ where $\a:A\ra B$ is an isogeny of abelian schemes,
$L$ is a relatively ample, symmetric,  line  bundle  on  $B$
trivialized along the zero section. Let $d=\rk\pi_*L$.
The main result of this section is the following theorem.
               
\begin{thm}\label{isogmain}
\begin{enumerate}
\item One has
$$gcd(12,\deg(\a))\cdot(\De(\a^*L)-\deg(\a)\cdot\De(L))=0.$$
\item If $\deg(\a)$ is odd then
$$\det\pi_*(\a^*L)=\deg(\a)\cdot\det\pi_*L
+d\cdot\det\pi_*\O_{\ker{\a}}+\zeta$$
where $gcd(3,\deg(\a))\cdot\zeta=0$.
\end{enumerate}
More precisely, these equalities in $\Pic(S)$
are realized by  canonical isomorphisms
of line bundles, compatible with arbitrary base changes.
\end{thm}
           
\begin{rem} In the case $d=1$ this follows from the result of
Moret-Bailly in \cite{MB}, VIII, 1.1.3. When $d$ is even the
second equality of the theorem can be rewritten as
$$gcd(3,\deg(\a))\cdot(\De'(\a^*L)-\deg(\a)\cdot\De'(L))=0$$
where $\De'(L):=\det\pi_*L+\frac{d}{2}\cdot\ov{\om}_A$.
\end{rem}
           
Outside of characteristic $3$ we can
improve Theorem \ref{isogmain} for some isogenies of degree
divisible by $3$.
           
\begin{thm}\label{isogmain2} Assume that $3$ is not among residue
characteristics of $S$ and that $d$ is relatively prime to $3$.
Assume that $3^k\cdot\ker(\a)=0$, that $K(\a^*L)^{(3)}$ is
annihilated by $\frac{1}{3}\cdot|K(\a^*L)^{(3)}|^{\frac{1}{2}}$,
and that $|K(\a^*L)^{(3)}/(\ker(\a)+3K(\a^*L)^{(3)})|>3$.
Then one has a canonical isomorphism
of line bundles on $S$ realizing the equality
$$\De(\a^*L)=\deg(\a)\cdot\De(L)$$
in $\Pic(S)$.
\end{thm}
         
\begin{rem} When $d$ is even the similar relation holds for
elements $\De'(\cdot)$ instead of $\De(\cdot)$.
\end{rem}
                                         
\begin{cor} Let $L$ be an ample symmetric line bundle on an
abelian scheme $A/S$, trivialized along the zero section.
Then for any odd integer $n>0$ which is not divisible by $3$
one has
$$\De(L^{n^2})=n^{2g}\cdot\De(L).$$
Also one has
$$\De(L^4)^{(3)}=4^{g}\cdot\De(L)^{(3)}$$
and
$$\De(L^9)^{(2)}=9^{g}\cdot\De(L)^{(2)}.$$
Furthermore, if $g>1$ and $3$ is prime to $d=\rk\pi_*L$ and
to the residue characteristics of $S$ then
$$\De(L^9)=9^{g}\cdot\De(L).$$
\end{cor}

We are going to use the relative Fourier-Mukai transform,
so let us briefly recall some of its properties. For details
the reader should consult \cite{Muk2} and \cite{La}.
The Fourier-Mukai transform is the functor
$$\FF_A=\FF_{A/S}:\D^b(A)\ra\D^b(\hat{A}):
X\mapsto p_2^*(p_1^*X\ot\PP),$$
where $\PP$ is the (normalized) relative Poincar\'e bundle
on $A\times_S\hat{A}$.
It is compatible with arbitrary base changes and satisfies the following
fundamental property:
$$\FF_{\hat{A}}\circ \FF_{A}\simeq
(-\id_A)^*(\cdot)\ot\pi^*\ov{\om}_A^{-1}[-g]$$
where $g$ is the relative dimension of $\pi$. Because of this
the relative canonical bundles of $A$ and $\hat{A}$ often appear
when working with the Fourier-Mukai transform so it is useful to know
that there is a canonical isomorphism
$\ov{\om}_{\hat{A}}\simeq\ov{\om}_A$
(see \cite{La}, 1.1.3).
              
Also for any homomorphism $f:A\ra B$ of abelian schemes over $S$ one
has the following canonical isomorphisms
\begin{equation}\label{hom1}
\FF_B\circ f_*\simeq \hat{f}^*\circ\FF_A,
\end{equation}
\begin{equation}\label{hom2}
\FF_A\circ f^!\simeq \hat{f}_*\circ\FF_B.
\end{equation}
The particular case of (\ref{hom1}) is the isomorphism
$$e^*\circ\FF_A\simeq\pi_*$$
where $e:S\ra A$ is the zero section, $\pi:A\ra S$ is the projection.

The  following  lemma  will  be used in the proof of Theorem
\ref{isogmain2}.

\begin{lem}\label{fur2}
For any relatively ample line bundle $L$ on $B$ trivialized along
the zero section one has a canonical isomorphism
$$\phi_L^*\FF_B(L)\simeq\pi^*\pi_*L\ot L^{-1}.$$
\end{lem}
            
\Pf. Making the base change $\phi_L:B\ra\hat{B}$ of the
projection $p_2:B\times\hat{B}\ra\hat{B}$ we can write
$$\phi_L^*\FF_B(L)\simeq\phi_L^*p_{2*}(p_1^*L\ot\PP)
\simeq p_{2*}(p_1^*L\ot (\id_B,\phi_L)^*\PP)$$
where in the latter expression $p_2$ denotes the projection of the product
$B\times_S B$ on the second factor. But we have an isomorphism
$$\mu^*L\simeq p_1^*L\ot p_2^*L\ot (\id_B,\phi_L)^*\PP$$
since a trivialization of $L$ along the zero section is
equivalent to a cube structure on it (see \cite{Br}). Hence,
$$\phi_L^*\FF_B(L)\simeq p_{2*}(p_2^*L^{-1}\ot\mu^*L)
\simeq L^{-1}\ot p_{2*}\mu^*(L)\simeq L^{-1}\ot\pi^*\pi_*L$$
as required.
\ed

\vspace{3mm}
   
\noindent
{\it Proof of Theorem \ref{isogmain}}.
             
Applying (\ref{hom2}) to the isogeny $\a:A\ra B$ we obtain
$$\FF_A(\a^!L)\simeq \hat{\a}_*\FF_B(L).$$
Restricting this isomorphism to the zero section we obtain
$$\pi_*(\a^!L)\simeq e^*\FF_A(\a^!L)\simeq
\pi_*(\FF_B(L)|_H)$$
where $H=\ker(\hat{\a})\sub\hat{B}$.
Taking the determinant of this isomorphism we get
$$\det\pi_*(\a^!L)=\det\pi_*(\FF_B(L)|_H).$$
Using the fact that $\a^!L\simeq\a^*L\ot\om_{A/S}\ot\om_{B/S}^{-1}$
we can rewrite the left hand side as follows:
$$\det\pi_*(\a^!L)=\det\pi_*(\a^*L)+n\cdot d\cdot\om_{\a}$$
where $\om_a=\ov{\om}_A-\ov{\om}_B$, $n=\deg\a$.
   
Recall (see e.~g. \cite{De}) that for any vector bundle $E$ of rank $r$ on
$H$ one has
\begin{equation}\label{detpi}
{\det}\pi_*(E)=r\cdot{\det}\pi_*\O_H +\Nm_{H/S}({\det} E)
\end{equation}
where $\Nm_{H/S}:\Pic(H)\ra\Pic(S)$ is the corresponding norm
homomorphism.
   
Consider the line bundle
$M=\det\FF_B(L)\ot\pi^*(\det\pi_*L)^{-1}$ on $\hat{B}$.
Then $M$ has a canonical trivialization along the zero section.
Furthermore, since the Fourier transform commutes with $[-1]^*$,
we have an isomorphism $[-1]_{\hat{B}}^*M\simeq M$ compatible with the
symmetry structure on $L$. Applying
the equation (\ref{detpi}) to $E=\FF_B(L)|_H$ we get
\begin{align*}
&\det\pi_*(\FF_B(L)|_H)=
d\cdot\det\pi_*\O_H+\Nm_{H/S}(\pi^*\det\pi_*L+M|_H)=\\
&d\cdot\det\pi^*\O_H+ n\cdot\det\pi_*L+\Nm_{H/S}(M|_H).
\end{align*}
Hence,
\begin{equation}\label{is1}
\det\pi_*(\a^*L) +n\cdot d\cdot\om_{\a}=
n\cdot\det\pi_*L+\Nm_{H/S}(M|_H)+d\cdot\det\pi_*\O_H.
\end{equation}
This implies that
$$\De(\a^*L)-n\cdot\De(L)=2\cdot\Nm_{H/S}(M|_H)+
2\cdot d\cdot\det\pi_*\O_H - n\cdot d\cdot\om_{\a}.$$
Recall that one has the canonical isomorphism
$$(\pi_*\O_H)^{\du}\simeq\om_H\ot\pi_*\O_H$$
where $\om_H=(\pi_*\O_{\hat{H}})^H$
is the line bundle of (relative) invariant measures on
$H$. Passing to determinants we get
$$2\det\pi_*\O_H=-n\cdot\om_H.$$
Hence,
$$\De(\a^*L)-n\cdot\De(L)=2\cdot\Nm_{H/S}(M|_H)-
n\cdot d\cdot\om_H - n\cdot d\cdot\om_{\a}=
2\cdot\Nm_{H/S}(M|_H)$$
since $\om_H\simeq\om_{\hat{\a}}\simeq\om_{\a}^{-1}$.
Also using that $\ker(\a)$ is Cartier dual to $H$ we deduce
from (\ref{is1}) the following equality
$$\det\pi_*(\a^*L)=n\cdot\det\pi_*L+
d\cdot\det\pi_*\O_{\ker\a}+\zeta$$
where $\zeta=\Nm_{H/S}(M|_H)$.
It remains to show that $gcd(24,2n)\cdot\Nm_{H/S}(M|_H)=0$
and that $gcd(3,n)\cdot\Nm_{H/S}(M|_H)=0$ if $\deg(\a)$ is odd.
To this end let us decompose $H$ into a product of
two group schemes $H\simeq H'\times_S H''$ such that
the order of $H'$ is odd, while the order of $H''$
is a power of 2. Then the cube structure on $M$ induces
the decomposition of $M|_H$ into the external
tensor product of $M|_{H'}$ and $M|_{H''}$. Hence, we obtain
$$\Nm_{H/S}(M|_H)=|H''|\cdot\Nm_{H'/S}(M|_{H'})+
|H'|\cdot\Nm_{H''/S}(M|_{H''}).$$
Recall  that  since   $M$   has   a   cube   structure   and
$[-1]^*M\simeq M$, it follows that $[n]^*M\simeq M^{n^2}$
for any $n$. Now the multiplication by $3$ is an automorphism of $H''$,
hence
$$\Nm_{H''/S}(M|_{H''})=\Nm_{H''/S}([3]_{\hat{B}}^*M|_{H''})=
\Nm_{H''/S}(M^{9}|_{H''})=9\cdot\Nm_{H''/S}(M|_{H''}).$$
Thus, $8\cdot\Nm_{H''/S}(M|_{H''})=0$.
Similarly, using the multiplication by $2$ we
obtain
$$\Nm_{H'/S}(M|_{H'})=\Nm_{H'/S}([2]_{\hat{B}}^*M|_{H'})=
\Nm_{H'/S}(M^{4}|_{H'})=4\cdot\Nm_{H'/S}(M|_{H'}),$$
hence $3\cdot\Nm_{H'/S}(M|_{H'})=0$.
It remains to note that
\begin{align*}
&2\cdot\Nm_{H/S}(M|_{H})=\Nm_{H/S}(M|_{H})+
\Nm_{H/S}([-1]_{\hat{B}}^*M|_{H})=\\
&\Nm_{H/S}((M\ot [-1]_{\hat{B}}^*M)|_{H})=
\Nm_{H/S}((\id_{\hat{B}},\phi_M)^*\PP|_{H})
\end{align*}
due to an isomorphism
$M\ot [-1]_{\hat{B}}^*M\simeq(\id_{\hat{B}},\phi_M)^*\PP$.
In particular, since $H$ is annihilated by $n$ it follows that
$$2n\cdot\Nm_{H/S}(M|_{H})=0.$$
\ed

\begin{lem}\label{norm}
Let $B$ be an abelian scheme over $S$,
where $3$ is prime to the residue characteristics of $S$,
$L$ be an ample line
bundle on $B$ trivialized along the zero section, $H\sub B_3$
be a finite flat subgroup. Assume that $H$ is isotropic with respect to
the symplectic form $e^{L^3}$ on $K(L^3)$ and that
$|H|>3$. Then
there is a canonical isomorphism of $\Nm_{H/S}(L|_H)$ with the
trivial line bundle on $S$.
\end{lem}
           
\Pf . Note that $L|_H$ is annihilated by $3$ in
$\Pic(H)$. Indeed, the cube structure on $L$
gives an isomorphism
\begin{equation}\label{triv0}
L_{3x}\simeq L_x^3 \ot\lan x,x\ran^3
\end{equation}
where $\lan\cdot,\cdot\ran=(\id_B,\phi_L)^*\PP$
is the symmetric biextension of $B\times B$
associated with $L$ (see \cite{Br}). When $3x=0$ this gives a trivialization
of $L_x^3$.
Let us consider the finite
\'etale covering  $c:S'\ra S$ corresponding to a choice of a non-trivial
point  $\si\in H$. The degree of this covering is prime to $3$,
so it suffices to prove the triviality
of $\Nm_{H/S}(L|_H)$ after making the corresponding base change.
Thus, we can assume that we have a non-trivial
$S$-point $\si:S\ra H$. To compute
$\Nm_{H/S}(L)$ we decompose the projection $H\ra S$ into
the composition $H\ra \ov{H}\ra S$
where $\ov{H}=H/\lan\si\ran$,
$\lan\si\ran\sub H$ is the cyclic subgroup in $H$ generated by $\si$.
Now we claim that
\begin{equation}\label{si}
\Nm_{H/\ov{H}}(L)\simeq \pi^*(\si^*L\ot(-\si)^*L)
\end{equation}
where $\pi$ is the projection to $S$. Indeed, to give an isomorphism
of line bundles on $\ov{H}$ is the same as to give an isomorphism of
their pull-backs to $H$ compatible with the action of
$\lan\si\ran\sub H$. Now the cube structure on $L$
gives the following canonical isomorphism
\begin{equation}\label{triv1}
L_x \ot L_{x+\si}\ot L_{x-\si}\simeq L_x^3\ot L_{\si}\ot L_{-\si}.
\end{equation}
Composing it with the trivialization of $L_x^3$ obtained above
we get an isomorphism
\begin{equation}\label{triv}
L_x \ot L_{x+\si}\ot L_{x-\si}\simeq L_{\si}\ot L_{-\si}.
\end{equation}
It remains to check that this isomorphism is compatible with
the action of $\Z/3\Z$
(the action on the right hand side being trivial),
i.~e. that we get the same isomorphism making the
cyclic permutation of the left hand side and applying (\ref{triv})
to $x+\si$. One can check that the only cause for these isomorphisms
to be different is the difference between the two trivializations of
$\lan x,\si\ran^3$: the one is obtained using that $3\cdot x=0$
and the other is obtained from $3\cdot\si=0$. But this difference
is equal to $e^{L^3}(x,\si)$ which we assumed to be trivial.
Here are more details of this computation.
First note that the isomorphism (\ref{triv1})  only  uses
the cube structure and the fact that $3\cdot\si=0$.
If we consider the similar isomorphism with $x$ replaced by
$x+\si$ it will be compatible with the natural isomorphism
$\a_{\si}:L_{x+\si}^3\simeq L_x^3$ which holds for any cube structure
and any $\si$ such that $3\cdot\si=0$.
Thus, we have to compute the difference  between  two
trivializations of $L_{x+\si}^3$: one which is  obtained  by
directly   applying   $(\ref{triv0})$  to  $x+\si$  and  then
trivializing $\lan x+\si,x+\si\ran^3$, and the other
which   is   the   composition   of   $\a_{\si}$   with  the
similar trivialization of  $L_x^3$.  Note  that  $\a_{\si}$  can  be
obtained using the isomorphisms $(\ref{triv0})$ and
the isomorphism
$$\a'_{\si}:\lan x+\si,x+\si\ran^3\simeq\lan x,x\ran^3$$
which again only depends on the fact that $3\cdot\si=0$.
Let us denote by $\b_x:\lan x,x\ran^3\ra 0$ the natural
trivialization for $x\in B_3$.
Now our claim follows from the equality
$\b_{x+\si}=e^{L^3}(x,\si)^{\pm 1}\cdot\b_x\circ\a'_{\si}.$
     
Thus, the isomorphism (\ref{triv}) is compatible with the action of
$\Z/3\Z$, hence the isomorphism (\ref{si}).
But $\si^*L\ot(-\si)^*L$  is annihilated by $3$ in
$\Pic(S)$ and the degree of the projection $\pi:\ov{H}\ra S$ is
divisible by $3$ (here we use the assumption that $|H|>3$). Thus,
we obtain
$$\Nm_{H/S}(L)=\Nm_{\ov{H}/S}\Nm_{H/\ov{H}}(L)=
\Nm_{\ov{H}/S}(\pi^*\si^*L^2)=0$$
as required.
\ed

\noindent
{\it Proof of Theorem \ref{isogmain2}}.
Let $H\sub B_{3^k}$ be the preimage of $\ker(\hat{\a})\sub\hat{B}_{3^k}$
under the isomorphism $\phi_L|_{B_{3^k}}:B_{3^k}\ra\hat{B}_{3^k}$.
As the proof of Theorem \ref{isogmain} shows we only have to check the
triviality of the norm of $M$ restricted to $\ker(\hat{\a})$.
Since $\phi_L^*M\simeq L^{-d}$ by Lemma \ref{fur2},
this is equivalent to proving the
triviality of $\Nm_{H/S}(L|_H)$. Consider the subgroup $K=\a^{-1}(H)\sub A$.
Then $K\sub A_{3^{2k}}$ and the definition of $H$
implies that $K=K(\a^*L)^{(3)}$. Also since $d$ is prime to $3$, it follows
that $\ker(\a)$ is a maximal isotropic subgroup in $K$.
Now  we  claim  that  after making an \'etale base change of
degree prime
to $3$ we can find a finite flat subgroup $K_1\sub K$
containing $\ker(\a)$, with the following two properties:
\begin{enumerate}
\item $K$ is annihilated by $3\cdot |K_1/\ker(\a)|$,
\item the quotient $K/K_1$ is annihilated by $3$ and
has order $>3$.
\end{enumerate}
Indeed, since $|K/\ker(\a)|=|K|^{\frac{1}{2}}$,
all we need is to find $K_1$ such that
$\ker(\a)+3K\sub K_1\sub K$ and $|K/K_1|=9$. To get
such $K_1$ we just make an \'etale covering of $S$
(of degree prime to $3$)
corresponding to a choice of a subgroup of index $9$
in $K/(\ker(\a)+3K)$.
          
Let us denote
$K_1'=K_1/\ker(\a)\sub H\sub B$ and $H'=K/K_1=H/K'_1$.
Consider the isogeny
$f:B\ra B'=B/K'_1$ and let $L'$ be a line bundle on $B'$ defined by
$$L'=\Nm_{B/B'}(L)\ot\pi^*\Nm_{K'_1/S}(L|_{K'_1})^{-1}.$$
Then $L'$ is trivialized along the zero section and we have
$$\Nm_{H/S}(L|_H)=\Nm_{H'/S}(\Nm_{B/B'}L|_{H'})=
\Nm_{H'/S}(L'|_{H'})+|H'|\cdot \Nm_{K'_1/S}(L|_{K'_1}).$$
The latter term is trivial, since $\Nm_{K'_1/S}(L|_{K'_1})$
is annihilated by $3$ (see the proof of Theorem \ref{isogmain}).
Thus, it remains to prove the triviality of $\Nm_{H'/S}(L'|_{H'})$.
Since $|H'|>3$ we can apply Lemma \ref{norm} to $L'$ and $H'$
provided that $H'$ is isotropic with respect to the
standard symplectic form on $K(L^{\prime 3})$. This is equivalent
to asking that $H$ is isotropic in $K((f^*L')^3)$. But
the symplectic structure on the latter group is determined by the
the polarization associated with $(f^*L')^3$ (see \cite{MuAb}).
Now since $f^*\Nm_{B/B'}(L)$ is algebraically equivalent to
$L^{\deg(f)}$ we obtain that
$K((f^*L')^3)=K(L^{3\deg(f)})$ as symplectic groups. Now
$H\sub K(L^{3\deg(f)})$ is isotropic iff
$K=\a^{-1}(H)\sub  K(\a^*L^{3\deg(f)})$  is  isotropic.  But
$K\sub K(\a^*L)$,
so this follows from the fact that $K$ is annihilated by
$3\deg(f)=3\cdot |K'_1|$ by assumption.
\ed
              
\section{The method of Faltings and Chai}
\label{meth}
                                           
In this section we start proving Theorem \ref{main1}.
                                                             
Fix an odd prime number $p$.
Following the proof of Faltings and Chai we
consider the homomorphism $f_p:\Z\ra\Pic(S)^{(p)}$, such that
$f_p(n)$ is the $p$-primary component of $\De(L^n)$.
This is a "polynomial" function in $n$, which means that $\de^if_p=0$ for
some $i$ where $\de$ is the difference operator:
$\de\phi(n)=\phi(n+1)-\phi(n)$. As was noticed in \cite{FC} this
can  be seen by embedding $A$ into the product of projective
bundles  $\P(\pi_*L^a)\times_S\P(\pi_*L^b)$  for  relatively
prime  $a$  and $b$ (see Lemma \ref{degree} below for a more
precise result). This implies immediately that
the image of $f_p$ belongs to some finitely generated subgroup of
$\Pic(S)^{(p)}$ and that $f_p(n+p^N)=f_p(n)$ for sufficiently large $N$.
By Serre duality one has
$$f_p(1)+(-1)^g\cdot f_p(-1)=0$$
where $g$ is the relative dimension of $A/S$.
Thus, if we find an integer
$k\equiv -1\mod(p^N)$ such that $f_p(k)=k^g\cdot f_p(1)$
this would imply that $\De(L)^{(p)}=0$. This is always possible when
$p\equiv  1\mod(4)$.
Indeed, we claim that if $p>3$ then
$$f_p(m^2n)=m^{2g}f_p(n)$$
for all $n$ and $m\neq 0$.
This follows immediately from Theorem \ref{isogmain} applied
to the isogeny $[m]_A:A\ra A$ and the line bundle $L^n$
(we don't have to worry about the factor $12$ since we only consider
$p$-primary component of the equality of Lemma \ref{degree} and $p>3$).
If $p\equiv 1\mod(4)$ then we can find $k=m^2\equiv -1\mod(4)$, so we are
done.
Hence,  we  can  assume  that $p\equiv
-1\mod(4)$. In this case one can always find some integers
$n$ and $m$ such that $n^2+m^2\equiv -1\mod(p^N)$.
Now let us consider the isogeny $\a:A^2\ra A^2$ given by the
matrix
$\left( \matrix {[n]_A} & {[m]_A} \\ {[-m]_A} & {[n]_A}
\endmatrix\right)$.
Then it is easy to see that $\a^*(L\boxtimes L)\simeq
L^k\boxtimes  L^k$ where $k=n^2+m^2$.
Note that possibly changing initial $n$
and $m$ we can achieve that $k$ is prime to $3$.
Applying Theorem \ref{isogmain} we find that
$$4\cdot\De(L^k\boxtimes L^k)=
4\cdot k^{2g}\De(L\boxtimes L).$$
Hence, $d\cdot f_p(k)=d\cdot k^g\cdot f_p(1)$. As we noticed
above this implies that $d\cdot\De(L)^{(p)}=0$.
This finishes  the  first  step in the proof of Theorem
\ref{main1}.
              
Now let us prove that $\De(L)^{(p)}=0$ for $p\ge 2g+1$, $p\neq 3$.
The only new ingredient we need is the following lemma.
Let us say that a function $\phi:\Z\ra G$, where $G$  is  an
abelian group, has degree $\le l$ if $\de^{l+1}\phi=0$
where $\de\phi(n)=\phi(n+1)-\phi(n)$.
              
\begin{lem}\label{degree}
Let $\pi:X\ra S$ be a smooth projective morphism
of pure dimension $g$, $L$ be a line bundle on $X$.
Then the function $f:\Z\ra\Pic(S)$ defined by
$f(n)=\det\pi_*(L^n)$ has degree $\le g+1$.
\end{lem}
              
\Pf . This follows from Elkik's construction (based on
ideas of Deligne in \cite{De}), see \cite{Elkik}, IV.1.3.
\ed

Applying this lemma to our abelian scheme $A/S$
and the line bundle $L$ on it
we deduce that $f_p$  has  degree  $\le g+1$.
Now the vanishing of $\De(L)^{(p)}$ for $p\ge 2g+1$, $p\neq 3$,
is implied by the following lemma.
              
\begin{lem}\label{ar1}
Let $p$ be a prime number, such that $p\ge 2g+1$,
$p\neq 3$, $\phi:\Z\ra\Z/p^k\Z$ be a function of degree
$\le g+1$, such that $\phi(n+p^N)=\phi(n)$ for sufficiently
large $N$. Assume that
$\phi(m^2n)=m^{2g}\cdot \phi(n)$ for all $n$ and
$m$. Then
$\phi(n)=n^g\cdot \phi(1)$ for all $n$.
\end{lem}
              
\Pf . Replacing $\phi$ by $\phi(n)-n^g\cdot \phi(1)$
we can assume that $\phi(1)=0$. In this case the assertion
of Lemma is that $\phi=0$. An easy induction in $k$ shows that
it suffices to prove this for $k=1$.
Then we can find a polynomial
$\phi'(x)\in\Z/p\Z[x]$ of degree $\le p-1$ such that
$\phi'(n)=\phi(n)$ for $n=0,1,\ldots,p-1$. Since $\phi$ is the
function of degree $\le p-1$, it is determined uniquely
by the set of its $p$ consequtive values. The same is true for
$\phi'$ considered as a function $\Z\ra\Z/p\Z$. It follows
that $\phi'(n)=\phi(n)$ for all $n$, in particular, $\phi(n)$
depends only on $n\mod(p)$.
Let us fix a non-quadratic residue
$a$ modulo $p$. We know that $\phi(n)=0$ if $n$ is a square
modulo $p$, and that
$\phi(n)=a^{-g}\cdot    n^g\cdot\phi(a)$   if   $n$   is
not a square modulo $p$. Hence, for some  $\la\in\Z/p\Z$  we
have
$$\phi(n)=\la\cdot n^g\cdot(1-n^{\frac{p-1}{2}})$$
for  all  $n$. Now if $\la\neq 0$ then
the right hand side is given by a polynomial of
degree $g+\frac{p-1}{2}\le p-1$. Therefore, we actually
have an identity of polynomials in $\Z/p\Z[x]$
which implies that $\deg(\phi)=g+\frac{p-1}{2}$. But
this contradicts to $\deg(\phi)\le g+1$, hence, $\la=0$
as required.
\ed
             
\begin{rem} The fact that the element
$\De(L)\in\Pic(S)$ has finite order is proved in \cite{FC}
along the same lines. One should consider the function
$f_0:\Z\ra\Pic(S)/\Pic(S)^{tors}$ where $\Pic(S)^{tors}$
is the torsion subgroup of $\Pic(S)$, such that
$f_0(n)=\De(L^n)\mod\Pic(S)^{tors}$. Then $f_0$ is a
function of finite degree, hence, its image is a finitely
generated free group. Then the identity
$f_0(n^2)=n^{2g}\cdot f_0(1)$ for infinitely many $n$ implies that
$f_0(n)=n^g\cdot f_0(1)$ for all $n$. Applying this
to $n=-1$ and using Serre duality as above we deduce that
$f_0=0$. At last, the bound on the 2-primary torsion of
$\De(L)$ is obtained by considering the isogeny $A^4\ra A^4$
given by a $4\times 4$ matrix of multiplication by a
quaternion $n+m\cdot i+p\cdot j+q\cdot k$ such that
$n^2+m^2+p^2+q^2\equiv -1(N)$ for sufficiently divisible $N$
(see \cite{FC}).
\end{rem}

\section{Some arithmetics}
\label{arsec}
            
Let $\phi:\Z\ra\Z/p^k\Z$ be a map, $g\ge 1$ be an integer. Let
us  say  that $\phi$ is $g$-{\it special}, if $\phi$ has degree
$\le g+1$ and $\phi(m^2n)=m^{2g}\cdot\phi(n)$ for all $n$ and $m$.
In particular, since $\phi$ has  finite  degree  it  factors
through $\Z/p^N\Z$ for sufficiently large $N$.
In this terminology Lemma \ref{ar1} says that
for $p>3$, $g\le\frac{p-1}{2}$ any $g$-special
map has form $\phi(n)=n^g\cdot \phi(1)$. In this section we'll
study $g$-special maps for other values of $g$.
Our main result is the following theorem, which combined
with  results  of  the  previous  section  implies the
first part of Theorem \ref{main1}, except for the fact
that $\De(L)^{(p)}$ is  annihilated  by  $p^2$  if  $p=2g-1$.
The latter statement will be proved together with the second
half of Theorem \ref{main1} in the end of this section.
            
\begin{thm}\label{ar2}
If $p>3$ and $g<\frac{3p-1}{2}$, $g\neq\frac{p+1}{2}$,
then for any
$g$-special map $\phi:\Z\ra\Z/p^k\Z$
one has $p\cdot \phi(n)=p\cdot n^g\cdot \phi(1)$.
\end{thm}
            
We will use the condition that the degree of $\phi$ is
$\le g+1$  in  the  following  form.  Let  us  consider  the
generating function
$$F(t)=\sum_{n\ge 0}\phi(n)t^n\in\Z/p^k\Z[[t]].$$
Then the condition $\deg(\phi)\le g+1$ implies that
$$F(t)\cdot (t-1)^{g+2}=P(t)$$
where $P(t)\in\Z/p^k\Z[t]$ is a polynomial in $t$.
In particular, if $\phi(n+p^N)=\phi(n)$ for all $n$
then
$$F(t)=Q_{\phi}(t)\cdot (1-t^{p^N})^{-1}$$
where $Q_{\phi}(t)=\sum_{n=0}^{p^N-1}\phi(n)$. Thus,
$Q_{\phi}(t)\cdot (t-1)^{g+2}$ is divisible by $t^{p^N}-1$
in $\Z/p^k\Z[t]$.
     
We are particularly interested in the case $k=1$.
In this case we obtain that $Q_{\phi}(t)$ is
divisible by $(t-1)^{p^N-g-2}$.
Let  us denote
$$S_r(t)=\sum_{n=0}^{p-1}n^r\cdot t^n\in\Z/p\Z[t]$$
for $r\ge 0$ (in case $r=0$ our convention is that $0^0=1$).
Note that for $r>0$ one has $S_r(t)=S_{r+p-1}(t)$.
For every polynomial $Q\in\Z/p\Z[t]$ we denote
by $v_{(t-1)}(Q)$ the maximal power of $(t-1)$ dividing $Q$.
            
\begin{lem}\label{val}
One has
$$v_{(t-1)}(S_r(t))=p-1-r$$
for $0\le r\le p-1$.
\end{lem}
             
\Pf . For $r=0$ we have an identity
$$S_0(t)=(t-1)^{p-1}$$
which follows from the congruence
${p-1 \choose i}\equiv (-1)^i\mod(p)$.
Now the identity $S_{r+1}(t)=t\cdot\frac{d}{dt}S_r(t)$
and an easy induction show that
$$S_r(t)\equiv (-1)^r\cdot r!\cdot t^r\cdot (t-1)^{p-1-r}
\mod((t-1)^{p-r})$$
for  $0\le r\le p-1$.
\ed
                       
The first step in the proof of Theorem \ref{ar2} is
the following lemma.
            
\begin{lem}\label{modp}
Let $\phi:\Z\ra\Z/p\Z$ be a $g$-special map.
If $p>3$ and $\frac{p-1}{2}<g<2p-1$ then
$$\phi(n)=\la\cdot n^g+\mu\cdot n^{g-\frac{p-1}{2}}\mod(p)$$
for some constants $\la,\mu\in\Z/p\Z$.
\end{lem}
            
\Pf. It is easy to see that for any $\la$ and $\mu$ the
map
$$n\mapsto \la\cdot n^g+\mu\cdot n^{g-\frac{p-1}{2}}\mod(p)$$
is $g$-special. Hence, if we write
$$\phi(n)=\la\cdot n^g+\mu\cdot n^{g-\frac{p-1}{2}}+\phi'(n)$$
for some $\la$ and $\mu$ then $\phi'$ will also be a
$g$-special map. Choosing $\la$ and $\mu$ appropriately
we can achieve that $\phi'(1)=\phi'(a)=0$ for
some $a$ which is a not a square modulo $p$. Replacing
$\phi$ by $\phi'$ we can assume that this condition holds
for $\phi$. Since $\phi(n+p)=\phi(n)$ for every
$n\not\equiv 0\mod(p)$ (this follows from $g$-speciality
and the fact that $(n+p)n^{-1}$ is a square in $\Z/p^N\Z$)
we deduce that $\phi(n)=0$ for all
$n\not\equiv 0\mod(p)$. On the other hand,
$\phi(p^2n)=p^{2g}\phi(n)=0$, hence the only non-trivial
values of $\phi$ are $\phi(pn)$ for $n\not\equiv 0\mod(p)$.
If $\deg\phi<p$ then this implies immediately that $\phi=0$,
so we can assume that $g\ge p-1$.
For all $n\not\equiv 0\mod(p)$ we have $\phi(pn+p^2)=(1+p\cdot n^{-1})^g
\cdot\phi(pn)=\phi(pn)$. In particular, $\phi$ depends only
on $p\mod(p^2)$ and
$$Q_{\phi}(t)=\sum_{n=0}^{p^2-1}\phi(n)t^n=
\sum_{n=1}^{p-1}\phi(pn)t^{pn}.$$
Now $n\mapsto\phi(pn)$ is a $p$-special map depending
only on $n\mod(p)$. Hence,
$$\phi(pn)=a\cdot n^g+b\cdot n^{g-\frac{p-1}{2}}$$
for some $a,b\in\Z/p\Z$. Therefore,
$$Q_{\phi}(t)=a\cdot S_g(t^p)+b\cdot S_{g-\frac{p-1}{2}}(t^p).$$
As we have seen above the fact that $\deg\phi\le g+1$
implies that $v_{(t-1)}(Q_{\phi})\ge p^2-g-2$.
Now we claim that for all $g$ such that $p-1\le g\le 2(p-1)$
the valuations of $S_g(t^p)$ and $S_{g-\frac{p-1}{2}}(t^p)$
at $(t-1)$ are less than $p^2-g-2$.
This would imply that $a=b=0$, hence $\phi=0$ as required.
To prove our claim let us apply Lemma \ref{val} to compute
$v_{(t-1)}S_g$ and $v_{(t-1)}S_{g-\frac{p-1}{2}}$.
For  $g=p-1$  we  get  $v_{(t-1)}S_g=0$, while for $p-1<g\le
2(p-1)$ we have $v_{(t-1)}S_g(t^p)=p\cdot(2(p-1)-g)<p^2-g-2$.
Similarly, for $p-1\le g\le\frac{3(p-1)}{2}$ we get
$v_{(t-1)}S_{g-\frac{p-1}{2}}(t^p)=
p\cdot(p-1-g+\frac{p-1}{2})<p^2-g-2$, while for
$\frac{3(p-1)}{2}<g\le 2(p-1)$ we have
$v_{(t-1)}S_{g-\frac{p-1}{2}}(t^p)=
p\cdot(2(p-1)-g+\frac{p-1}{2})$. Thus, to finish the proof
we need the inequality
$$p\cdot(\frac{5(p-1)}{2}-g)<p^2-g-2$$
in   this   case,   but   when   $p>3$   it   follows   from
$g>\frac{3(p-1)}{2}$.
\ed
            
\begin{rem}  In  fact,  for  $p>3$  one  can  prove that the
conclusion of  the  previous  lemma  remains  true when
$g$ belongs to one of the following intervals of integers:
$[2p,\frac{5}{2}(p-1)]$, $[\frac{5p+1}{2},3(p-1)]$,
$[3p,\frac{7}{2}(p-1)]$, $[\frac{7p+1}{2},4(p-1)]$, etc.
(for given $p$ the set of such $g$ is finite).
\end{rem}
                                    
Next step is to consider $g$-special maps depending
only on $n\mod(p^2)$.
            
\begin{lem}\label{modp2}
Let $\phi:\Z\ra\Z/p^2\Z$ be a $g$-special map such that
$\phi(n+p^2)=\phi(n)$ for all $n$. Assume that
$p>3$ and $\frac{p+1}{2}<g<\frac{3p-1}{2}$. Then
$$p\cdot\phi(n)=p\cdot n^g\cdot\phi(1)$$
for all $n$.
\end{lem}
             
\Pf  . Replacing $\phi$ by $\phi-n^g\cdot\phi(1)$ we
can assume that $\phi(1)=0$. In this case we need to show
that $\phi=0\mod(p)$. Applying Lemma \ref{modp} to $\phi\mod(p)$
we obtain
\begin{equation}\label{phipsi}
\phi(n)=c\cdot (n^g- n^{g-\frac{p-1}{2}})+
p\cdot\psi(n)
\end{equation}
for some constant $c\in\Z/p^2\Z$ and some map
$\psi:\Z\ra\Z/p\Z$.
Now for every $n\not\equiv 0\mod(p)$ we have
\begin{align*}
&\phi(n+p)-\phi(n)=(1+p\cdot n^{-1})^g\cdot\phi(n)-\phi(n)=
p\cdot g\cdot n^{-1}\cdot\phi(n)=\\
&p\cdot c\cdot g\cdot (n^{g-1}-n^{g-\frac{p+1}{2}}).
\end{align*}
On the other hand, subtracting the equation (\ref{phipsi})
for $n$ from that for $n+p$ we get
$$\phi(n+p)-\phi(n)=
p\cdot c\cdot(g\cdot n^{g-1}-
(g-\frac{p-1}{2})\cdot n^{g-\frac{p+1}{2}})+
p\cdot(\psi(n+p)-\psi(n)).$$
It follows that
$$\psi(n+p)-\psi(n)=-c\cdot\frac{p-1}{2}\cdot
n^{g-\frac{p+1}{2}}=\frac{c}{2}\cdot
n^{g-\frac{p+1}{2}}.$$
Hence, for every $m$ one has
\begin{equation}\label{psi1}
\psi(n+p\cdot m)=\psi(n)+\frac{c}{2}\cdot m\cdot
n^{g-\frac{p+1}{2}}
\end{equation}
provided that $n\not\equiv 0\mod(p)$.
            
Also we claim that
\begin{equation}\label{psi2}
\psi(p\cdot n)=\la\cdot n^g+\mu\cdot n^{g-\frac{p-1}{2}}
\end{equation}
for some $\la,\mu\in\Z/p\Z$.
Indeed, since $g>\frac{p+1}{2}$ the equation (\ref{phipsi})
shows that $\phi(p\cdot n)=p\cdot\psi(p\cdot n)$. Hence, the map
$n\mapsto\psi(p\cdot n)$ is $g$-special and the assertion follows
from Lemma \ref{modp}.
            
Now (\ref{phipsi}) shows that $\psi$ depends only on
$n\mod(p^2)$ and has degree $\le g+1$. Hence, the corresponding
polynomial $Q_{\psi}(t)=\sum_{n=0}^{p^2-1}\psi(n)t^n$
is divisible by $(t-1)^{p^2-g-2}$. Using (\ref{psi1}) and
(\ref{psi2}) we can write
\begin{align*}
&Q_{\psi}(t)=
\sum_{n=1}^{p-1}\sum_{m=0}^{p-1}\psi(n+p\cdot m)t^{n+p\cdot m}
+\sum_{n=1}^{p-1}\psi(p\cdot n)=\\
&(\sum_{n=1}^{p-1}\psi(n)t^n)\cdot S_0(t^p)+
\frac{c}{2}\cdot S_{g-\frac{p+1}{2}}(t)\cdot S_1(t^p)
+\la\cdot S_g(t^p)+\mu\cdot S_{g-\frac{p-1}{2}}(t^p).
\end{align*}
In the case $g\le p-2$ we have $p^2-g-2\ge p^2-p$,
hence, $v_{(t-1)}Q_{\psi}\ge p^2-p$. To prove that $c=0$
in this case it is sufficient to check that
$v(g):=v_{(t-1)}S_{g-\frac{p+1}{2}}(t)+v_{(t-1)}S_1(t^p)<p^2-p$
and that $v(g)$ differs from $v_{(t-1)}S_g(t^p)$ and
$v_{(t-1)}S_{g-\frac{p-1}{2}}(t^p)$. One can check using Lemma
\ref{val} that this is indeed the case.
When $g\ge p-1$ we can omit the  first  term  in  the  above
expression for $Q_{\psi}(t)$ when considering
$Q_{\psi}(t)\mod(t-1)^{p^2-g-2}$. Hence, to deduce that $c=0$
one should check using Lemma \ref{val} that
$v(g)<p^2-g-2$ and $v(g)$ differs from the valuations of two
other terms. We omit the details of this simple computation.
\ed
               
\noindent
{\it Proof of Theorem \ref{ar2}}.
The  case  $g\le\frac{p-1}{2}$  follows  from  Lemma
\ref{ar1} so we only consider $g>\frac{p+1}{2}$.
Also as usual we can assume that $\phi(1)=0$.
An easy induction in $k$ shows that it suffices to consider
the case $k=2$. In the latter case we have
$\phi(p^2\cdot n)=0$ for all $n$ and $\phi(n+p^2)=\phi(n)$ for
$n\not\equiv 0\mod(p)$.
According to Lemma \ref{modp} we can write
$$\phi(n)=c\cdot (n^g- n^{g-\frac{p-1}{2}})+
p\cdot\psi(n)$$
for some constant $c\in\Z/p^2\Z$ and some map
$\psi:\Z\ra\Z/p\Z$.
In particular, $\phi(p\cdot n)=p\cdot\psi(p\cdot n)$
and $n\mapsto\psi(p\cdot n)$ is a $g$-special
map. Now Lemma \ref{modp} implies that
$\psi(p\cdot(n+p))=\psi(p\cdot n)$, hence
$\phi(p\cdot n+p^2)=\phi(p\cdot n)$. Therefore,
$\phi(n)$ depends only on $n\mod(p^2)$ and we can
apply Lemma \ref{modp2} to finish the proof.
\ed
                         
Now we are turning to the proof of the second half of
Theorem \ref{main1}.
Note that for any prime $p$ the function
$f_p:n\mapsto\De(L^n)^{(p)}$ satisfies the
following property:
$$f_p(m^2)=m^{2g}\cdot f_p(1)$$
for all $m$ such that $m\not\equiv 0\mod(p)$. Indeed, changing
$m$ by $m+p^N$ if necessary we may assume that $m$ is odd and
is   prime   to   $3$,   hence  this  follows  from  Theorem
\ref{isogmain}.
          
\begin{prop} Let
$f:\Z\ra\Z/p^k\Z$ be a function of degree $\le g+1$, such
that $f(m^2)=m^{2g}\cdot f(1)$ for all $m$ such that
$m\not\equiv 0\mod(p)$. Then there exists an integer
$n(p,g)$ depending only on $p$ and $g$ such that
$p^{n(p,g)}\cdot (f(n)-n^g\cdot f(1))=0$ for
all $n$.
\end{prop}
   
\Pf . First of all, replacing $f$ by $(g+1)!\cdot f$
we can assume that
$$f(n)-n^g\cdot f(1)=a_0+n\cdot a_1+\ldots +n^{g+1}\cdot a_{g+1}$$
for some $a_i\in \Z/p^k\Z$. Let $n_0,\ldots,n_{g+1}$ be the
first $g+2$ positive integers of the form $m^2$ with
$m\not\equiv 0\mod(p)$. Then $f(n_i)-n_i^g\cdot f(1)=0$
by assumption, hence every coefficient $a_i$ is annihilated
by the Vandermonde determinant
$\De(n_0,\ldots,n_{g+1})=\prod_{i<j}(n_j-n_i)$.
It follows that we can take $n(p,g)$ to be
$v_p((g+1)!\cdot\De(n_0,\ldots,n_{g+1}))$.
\ed
   
The proof of the second part of Theorem \ref{main1} follows
immediately from this proposition: we can take
$N(g)=\prod_{p\le 2g-1}p^{n(p,g)}$.
To complete the proof of Theorem \ref{main1} it remains to prove
the statement concerning the prime  $p=2g-1$.  This  is  the
content of the following lemma.
   
\begin{lem} Let $p>3$ be a prime, $f:\Z\ra\Z/p^k\Z$ be
a $g$-special map where $g=\frac{p+1}{2}$. Then
$p^2\cdot (f(n)-n^g\cdot f(1))=0$.
\end{lem}
   
\Pf . Since $g+1<p$ we have
$$f(n)-n^g\cdot f(1)=a_0+n\cdot a_1+\ldots +n^{g+1}\cdot a_{g+1}$$
for some $a_i\in \Z/p^k\Z$. Let $n_0=0$,
$n_1,\ldots,n_{g+1}$ be the
first $g+1$ positive integers of the form $m^2$ with
$m\not\equiv 0\mod(p)$. Then $f(n_i)-n_i^g\cdot f(1)=0$
for all $i=0,\ldots,g+1$ by $g$-speciality,
therefore every coefficient $a_i$ is annihilated
by the Vandermonde determinant
$\De=\De(n_0,\ldots,n_{g+1})=\prod_{i<j}(n_j-n_i)$.
Since $g+1=\frac{p-1}{2}+2$, it follows that
$v_p(\De)=2$, hence the assertion.
\ed
   
\section{Eliminating of primes not dividing characteristic}
\label{charsec}
              
In the case when an odd prime $p$ is not among residue characteristics
of the base we can evaluate $\De(L)^{(p)}$
over $S[\frac{1}{p}]$ using the following result.
              
\begin{thm}\label{bc}
Let $L$ be a symmetric, relatively ample line
bundle over an abelian scheme $A/S$  trivialized  along  the
zero section. Let $p$ be an odd prime number which is  not  equal
to any of residue characteristics of $S$.
Then there exists a finite flat base change $c:S'\ra S$ of
degree prime to $p$, an isogeny of abelian $S'$-schemes
$\a:A'\ra B$, where $A'$ is obtained from $A$ by this base
change, such that $\deg(\a)$ is the power of $p$,
and a symmetric line bundle $M$ on $B$ together
with a symmetric isomorphism $\a^*M\simeq L'$ (where $L'$
is obtained from $L$ by the base change), such that
$\deg\phi_M$ is prime to $p$.
\end{thm}
              
Let us deduce Theorem \ref{main2} from this.
Choose a base change $c:S'\ra S$ as in  Theorem \ref{bc}.
Since $\Pic(S)^{(p)}\ra\Pic(S')^{(p)}$ is injective
and the construction of $\det\pi_*$ commutes with this base change we
can work with $A'$ instead of $A$. It remains to apply
Theorem \ref{isogmain} to isogeny $\a$ and Theorem \ref{main1}
to $M$ to deduce that $\De(L')^{(p)}=0$ if $p\neq 3$.
In the case $p=3$ by the same argument we always have
$3\cdot\De(L')^{(3)}=0$.  Now  if   the   $3$-type   of   the
polarization is different from $(1,\ldots,1,3^k)$, then
one can see easily from the construction of the isogeny $\a$ below
that the conditions of Theorem \ref{isogmain2} are satisfied
for $\a$ and $M$.
Hence, the triviality of $\De(L')^{(3)}$ in this case.

\vspace{3mm}
         
\noindent
{\it Proof of Theorem \ref{bc}}. We can assume that the base $S$
is connected. Let $K=K(L)^{(p)}$ be the $p$-primary component of the
finite flat group scheme $K(L)$ over $S$. Then
$K$ is \'etale over $S$ (since $p$ is not among residue characteristics
of $S$) and $\phi_L$ induces a skew-symmetric isomorphism $K\simeq \hat{K}$.
The fiber of $K$ over a geometric point of $S$ is a discrete
group of the form
$(\Z/p^{n_1}\Z)^2\times\ldots\times(\Z/p^{n_k}\Z)^2$
where the factors $(\Z/p^{n_i}\Z)^2$ are orthogonal to each other
with respect to a symplectic form,
$\Z/p^{n_i}\Z\times \{0\}\sub(\Z/p^{n_i}\Z)^2$
is a lagrangian subgroup for every $i$. Since $S$ is connected the
collection $n_1,\ldots,n_k$ doesn't depend on a point.
Let $n$ be the maximum of $n_1,\ldots, n_k$, so that $n$ is the
minimal number such that $p^n\cdot K=0$. Now let us construct
a canonical isotropic subgroup $I_0\sub K$, \'etale over $S$, such that
$K_0=I_0^{\perp}/I_0$ is annihilated by $p$. If $n=1$ we can take
$I_0=0$ so let's assume that $n\ge 2$. Then $p^{n-1}\cdot K$
is an isotropic subgroup in $K$ so we can consider the reduction
$\ov{K}=(p^{n-1}\cdot K)^{\perp}/(p^{n-1}\cdot K)$ with its induced symplectic
form. By induction we may assume that we already found an isotropic
subgroup $\ov{I}_0\sub K'$ such that $(\ov{I}_0)^{\perp}/\ov{I}_0$
is a $p$-group.
Now take $I_0$ to be the preimage of $\ov{I}_0$ in $K$.
              
Our base change $c:S'\ra S$ will be the finite flat covering corresponding
to a choice of a lagrangian subgroup in $K_0$. One
can construct such a covering in the following way. Let $p^{2r}$ be the
order of $K^0$. Start with a subscheme $\wt{S}$ in
$K_0\times_S\ldots\times_S K_0$ ($2r$ times) corresponding to symplectic
bases in $K_0$. Then $\wt{S}\ra S$ is a $\Sp_{2r}(\Z/p\Z)$-torsor.
Now let $P\sub\Sp_{2r}(\Z/p\Z)$ be the subgroup preserving the
standard $r$-dimensional lagrangian subgroup in $(\Z/p\Z)^{2r}$,
then we can take $S'=P\backslash \wt{S}$.
The degree of the covering $S'\ra S$ is equal to the number of
lagrangian subgroups in $(\Z/p\Z)^{2r}$ which is easily seen to be equal
to $\prod_{i=1}^r(p^i+1)$ (first compute the number of isotropic flags
and then divide by the number of flags in $(\Z/p\Z)^r$), which is prime
to $p$.
              
Let $(A',L',K',I'_0,K'_0)$ be the data obtained from
$(A,L,K,I_0,K_0)$ by the base change $S'\ra S$. Then
by construction we have a lagrangian subgroup $\ov{I}\sub K'_0$.
Taking its preimage by the morphism $(I'_0)^{\perp}/I'_0\ra K'_0$
we obtain a lagrangian subgroup $I\sub K'$. It remains to prove
that $L'$ descends to a symmetric bundle on $B=A'/I$.
Let $p:G\ra K'$ be the restriction of the Mumford's group of $L'$.
It is well-known that a descent of $L'$ to a line bundle on $A'/I$,
is equivalent to choosing a splitting of $p$ over $I$.
Thus, we are reduced to finding a trivialization of the group
extension $\G_m\ra p^{-1}(I)\ra I$, which is compatible with the isomorphism
$\tau:G\wt{\ra}(-\id)^*G$.
Note that such a trivialization is necessarily
unique since a homomorphism $f:I\ra\G_m$ satsifying $f(-x)=f(x)$
is trivial (recall that the order of $I$ is odd). Hence, it is sufficient
to prove local existence of such trivialization. Now locally
there exists a homomorphism $\si:I\ra p^{-1}(I)$ splitting $p$. Then
$x\mapsto \tau^{-1}\si(-x)$ is another such splitting, so that
$\tau^{-1}\si(-x)=\psi(x)\cdot \si(x)$ for some homomorphism
$\psi:I\ra\G_m$. Since the multiplication by $2$ is invertible
on $I$ there exists a homomorphism $\phi:I\ra\G_m$ such that
$\psi(x)=\phi(x)^2$. Then $x\mapsto \phi(x)\si(x)$ gives the
symmetric splitting.
\ed
             
\section{Complements}\label{comp}

\subsection{Case of elliptic curves}
Let us evaluate determinant bundles in the case $g=1$.
It is known (see e.g. \cite{Mu}) that in characteristics
$\neq 2,3$ the Picard group of the moduli stack of elliptic
curves (with one fixed point) is $\Z/12\Z$ and as generator
one can take the line bundle $\ov{\om}$ on this moduli stack
that associates to every family of elliptic curves $\pi:E\ra S$
the relative canonical bundle $\ov{\om}_E\in\Pic(S)$.
If $S$ is connected then
any symmetric line bundle $L$ on $E$, trivialized
along the zero section, is either isomorphic to
$L_d(e):=\O(d\cdot e)\ot\om_{E/S}^d$ where $e:S\ra E$ is the zero section,
or to $L_d(\eta):=\O((d-1)\cdot e+\eta)\ot\om_{E/S}^{d-1}$ where
$\eta:S\ra E$ is an everywhere non-trivial point of order 2.
Note that these line bundles are trivialized along the zero section,
since $e^*\O(e)\simeq\ov{\om}_E^{-1}$.
             
\begin{prop}\label{ell}
One has
$$\det\pi_*(L_d(e))=(\frac{d(d-1)}{2}+1)\cdot\ov{\om}_E,$$
$$\det\pi_*(L_d(\eta))=\frac{d(d-1)}{2}\cdot\ov{\om}_E.$$
In particular,
$$\De(L_d(e))=(d^2+2)\cdot\ov{\om}_E,$$
$$\De(L_d(\eta))=d^2\cdot\ov{\om}_E.$$
Furthermore, these equalities are represented by
canonical isomorphisms of line bundles.
\end{prop}
             
\Pf. Considering the push-forward of the exact sequence
$$0\ra\O((d-1)\cdot e)\ra\O(d\cdot e)\ra e_*e^*\O(d\cdot e)\ra
0$$
we deduce that
$$\det\pi_*\O(d\cdot e)-\det\pi_*\O((d-1)\cdot e)=
-d\cdot\ov{\om}_E.$$
Since $\pi_*\O(e)\simeq\O_S$ it follows that
$$\det\pi_*\O(d\cdot e)=(1-\frac{d(d+1)}{2})\ov{\om}_E,$$
hence $\det\pi_*L_d(e)=(\frac{d(d-1)}{2}+1)\ov{\om}_E$.
The case of $\O_d(\eta)$ is considered similarly using
the exact sequence
$$0\ra\O((d-1)\cdot e+\eta)\ra\O(d\cdot e+\eta)\ra
e_*e^*\O(d\cdot e+\eta)\ra 0$$
and the triviality of $\pi_*\O(\eta)$.
\ed
             
Note that $\De(L_d(e))$ gives a line bundle on the moduli stack
of elliptic curves $\AA_1$, while $\De(L_d(\eta))$ lives on the
stack $\wt{\AA}_1$ classifying elliptic curves
with a non-trivial point of order 2.
Now Proposition \ref{ell} combined with Theorem \ref{main1}
implies immediately that the order of $\ov{\om}$ in
$\Pic(\AA_1)$ is 12, while the order of the pull-back of
$\ov{\om}$ to $\wt{\AA}_1$ is 4. In particular,
$\De(L_3(e))=-\ov{\om}$ in $\Pic(\AA_1)$.
          
\subsection{Linear relations between
determinant bundles}
Let us first consider the determinant bundles $\det\pi_*(L^n)$
on  an  abelian  scheme $A/S$ of relative dimension $g=2$ or
$g=3$, where $L$ is a relatively
ample, symmetric line bundle on $A$ trivialized along the
zero section.
          
\begin{prop}\label{Ln}
If $g=2$ then one has
$$\det\pi_*L^n=\frac{4n-n^3}{3}\cdot\det\pi_*L+
\frac{n^3-n}{6}\cdot\det\pi_*L^2+
\frac{n(n-1)(n-2)}{6}\cdot d\cdot\ov{\om}_A,$$
where $d=\rk\pi_*L$. In particular,
$$\De(L^n)=\frac{4n-n^3}{3}\cdot\De(L)+\frac{n^3-n}{6}\cdot\De(L^2).$$
If $g=3$ then one has
$$\det\pi_*L^n=\frac{4n^2-n^4}{3}\cdot\det\pi_*L+
\frac{n^4-n^2}{12}\cdot\det\pi_*L^2+
\frac{n^2(n-1)(n-2)}{6}\cdot d\cdot\ov{\om}_A.$$
In particular, in this case
$$\De(L^n)=\frac{4n^2-n^4}{3}\cdot\De(L)+\frac{n^4-n^2}{12}\cdot
\De(L^2).$$
\end{prop}
         
\Pf . This is easily deduced from the fact that the function
$f:n\mapsto\det\pi_*L^n$ has  degree  $\le  g+1$.  Indeed,  by
Serre duality the values of this function at $n=-2,-1$
are expressed via those for $n=1,2$. Also $\det\pi_*\O_A=0$,
hence,  we  know  values  of $f$ at $n\in [-2,2]$ and we can
interpolate the rest.
\ed
         
\begin{cor} If $g=2$ or $g=3$, then
$$8\cdot 9\cdot\De(L)=4\cdot 9\cdot\De(L^2)=0.$$
\end{cor}
    
\Pf . This is proved by considering separately $2$-primary and
$3$-primary  parts  of  $\De(L)$  and  $\De(L^2)$, using the
previous proposition and Theorem \ref{isogmain}.
\ed
         
If $g=2$, $gcd(d,3)=1$ and $n$ is odd then we also get from
Proposition \ref{Ln} that $\De(L^n)=n\cdot\De(L)$.
If $g=3$, $gcd(d,3)=1$, and the characteristic is zero then according to
Kouvidakis one has $\De(L^2)=0$, hence in this case for
odd $n$ we get $\De(L^n)=\De(L)$. It would be interesting
to find similar dependences between $\De(L^n)$ in higher
dimensions  (see  section \ref{tor} for the case of even principal
polarization).
Recall (see Lemma \ref{degree}) that $n\mapsto\det\pi_*(L^n)$ is a function
of degree $\le g+1$. Hence, using Serre's duality,
one can express all $\det\pi_*(L^n)$
as   linear   combinations of $\ov{\om}$ and $\det\pi_*(L^i)$ where
$0\le i\le \frac{g}{2}+1$. However, we expect much more relations
between $\De(L^n)$. Here are some examples.
         
\begin{enumerate}
\item For $g\ge 2$ one has $\det\pi_*\O_A=0$.
\item Let $p$ be a prime, $p\equiv -1\mod(4)$. Then
for any $n$ such that $(n,p)=1$ one has
$$\De(L^n)^{(p)}=\left(\frac{n}{p}\right)\cdot
n^g\cdot\De(L)^{(p)}$$
where $\left(\frac{n}{p}\right)=\pm 1$ is the Legendre symbol.
\item Assume that $g\ge 2$.
Let $p$ be a prime such that $p\equiv -1\mod(4)$ and $p\ge (g+3)/2$.
Then one has $\De(L^n)^{(p)}\in\Z\De(L)^{(p)}$ for all $n$.
\item For odd $n$
one has $\De(L^{n+8})^{(2)}=((n+8)/n)^g\cdot\De(L^n)^{(2)}$.
In particular, if $n$ and $d$ are odd then
$\De(L^{n+8})^{(2)}=\De(L^n)^{(2)}$.
\end{enumerate}
         
(1) follows from the fact that
$R^i\pi_*\O_A=\We^i R^1\pi_*\O_A$. For the proof of (2) note
that for $(n,p)=1$ Theorem \ref{isogmain} implies that
$\De(L^n)^{(p)}=n^g\cdot\De(L)^{(p)}$ if $n$ is a square modulo $p$,
and that $\De(L^n)^{(p)}=(-n)^g\cdot\De(L^{-1})^{(p)}$ if $-n$ is a
square modulo $p$. But Serre's duality implies that
$\De(L^{-1})=-(-1)^g\cdot\De(L)$, hence  the  assertion.  To
prove (3) note that for $p\ge (g+3)/2$ all the elements $\De(L^n)^{(p)}$
are linear combinations of $\De(L^i)^{(p)}$ with $|i|<p$. It
remains to apply (1) and (2). At last, (4) follows
from Theorem \ref{isogmain}
since  $\frac{n+8}{n}$ is a square modulo resp. $2^k$.
                                    
\begin{prop} $d\cdot (\De(L^3)^{(2)}+3^g\cdot\De(L)^{(2)})=0$.
\end{prop}

\Pf . Let us denote
$\De(L,n)=\det\pi_*(L^n)-n^g\cdot\det\pi_*L$.
Then it is easy to see that
$$2\De(L,n)=\De(L^n)-n^g\cdot\De(L),$$
in particular, $\De(L,n)$ is a torsion element in $\Pic(S)$.
Note also that $n\mapsto\De(L,n)$ is a polynomial function.
Let us choose a sufficiently divisible integer $N>0$ such that
both  functions  $\De(L^n)$  and  $\De(L,n)$  of   $n$   are
$N$-periodic and are annihilated  by  $N$,  and  $N$  is
divisible by $6$. Now let $l$ be a prime, such that $N$ is
not divisible by $l$ and such that $l\equiv 3\mod(2^m)$
where $m>>0$. Then there exists a solution $(a,b)$ to the congruence
$a^2+lb^2\equiv -1\mod(N)$. Consider the isogeny
$\a:A^2\ra A^2$ given by the matrix
$\left( \matrix a & -lb \\ b & a \endmatrix\right)$.
Then it is easy to see that
$\a^*(L\boxtimes L^l)\simeq L^k\boxtimes L^{lk}$
where $k=a^2+lb^2$. In particular,
$\deg(\a)=k^{2g}$ and applying Theorem \ref{isogmain}
we obtain the equality
$$\det\pi_*(L^k\boxtimes L^{lk})=k^{2g}\cdot
\det\pi_*(L\boxtimes L^l).$$
Now using the equalities $\De(L,k)=\De(L,-1)$,
$\De(L,lk)=\De(L,-l)$, Serre's duality, the fact
that $k\equiv -1\mod(N)$, and the conditions on $N$ one
can easily deduce that
$$d\cdot(\De(L^l)+l^g\cdot\De(L))=0.$$
It remains to take the 2-primary part of this equality and
replace $l$ by $3$ in the obtained identity (this is
justified by the congruence $l\equiv 3\mod(2^m)$ with
$m>>0$).
\ed

\begin{cor} Assume that $n$ and $d$ are odd. Then
$$\De(L^n)^{(2)}=n^{g-1}\cdot\De(L)^{(2)}.$$
In particular, if $d=1$ and $g\ge 2$ then
$$\De(L^n)=n^{g-1}\cdot\De(L)$$
for any odd $n$.
\end{cor}

\Pf . Since $\De(L^n)^{(2)}$ for odd $n$ depends only
on $n\mod(8)$ the first equality follows from the case $n=3$ considered
above, Serre's duality, and the vanishing of
$4\cdot\De(L^n)^{(2)}$ for odd $n$. Now the second statement follows
from  the  fact  that  for  $d=1$  and  $g\ge  2$  one   has
$4\cdot\De(L^n)=0$ for any $n$.
\ed
 
\subsection{Torsion in the Picard group of moduli}\label{tor}
Let $\wt{\AA}_g$ be the moduli stack of the data
$(A/S,\Th)$ where $A/S$ is an abelian scheme of relative dimension $g$,
$\Th\sub A$ is an effective (relative) divisor which is symmetric and
defines a principal polarization of $A$.
One can normalize the line bundle $\O(\Th)$ over the universal abelian scheme
over $\wt{\AA}_g$ to obtain the line bundle $L$ which is trivial
along the zero section. In particular, we have an element
$\De(L)\in\Pic(\wt{\AA}_g)$.
Let $\wt{\AA}_g^+$ be
the irreducible component of $\wt{\AA}_g$ corresponding to even theta divisors,
$\wt{\AA}_g^+[\frac{1}{2}]$ be the localization of this stack over
$\Spec(\Z[\frac{1}{2}])$.
        
\begin{thm}\label{torsion} Assume that $g\ge 3$.
Then the torsion subgroup in $\Pic(\wt{\AA}_g^+[\frac{1}{2}])$
is isomorphic to $\Z/4\Z$ and is generated by $\De(L)$.
\end{thm}
        
\Pf . Since $\wt{\AA}_g^+$ has smooth geometrically
irreducible fibers over
$\Spec(\Z[\frac{1}{2}])$ (cf. \cite{FC}, IV 7.1)
it is sufficient to prove this statement in
characteristic zero.
Indeed, it is known that the order of $\De(L)$ is precisely 4 (cf.
\cite{MB2}), hence it would follow that $\De(L)$
generates the entire torsion subgroup in
the Picard group of the general fiber of $\wt{\AA}_g^+$.
Now it remains to prove that if some
line bundle over $\wt{\AA}_g^+$ is trivial over the general fiber then
it is trivial everywhere. If $\wt{\AA}_g^+$ were represented
by a scheme then we could apply
the argument from \cite{Mumst}, p. 103, to prove this. Since
it is not, we have to replace  $\wt{\AA}_g^+$  by a
$\PGL_N$-torsor over it which is representable
(see remark 1 after Theorem \ref{main2}),
apply the cited argument, and use the
fact that $\PGL_N$ has no non-trivial characters.
        
The corresponding analytic stack is the quotient (in the sense
of stacks) of the Siegel's half-space $\SH_g$ by the subgroup
$\Ga_{1,2}\sub\Sp_{2g}(\Z)$ consisting of matrices whose reduction modulo
2 preserves the standard even quadratic form $\sum_{i=1}^g x_iy_i$.
(cf. \cite{MB}, VIII, 3.4). It follows that
the torsion in the Picard group of this stack is an abelian group dual
to $\Ga_{1,2}/[\Ga_{1,2},\Ga_{1,2}]$ (cf. \cite{Mu}).
It remains to prove that the latter group is isomorphic to $\Z/4\Z$.
         
As is shown in \cite{Theta3}, Prop. 8.10, there is a normal subgroup
$\De\sub\Ga_{1,2}$ such that $\Ga_{1,2}/\De\simeq\Z/4\Z$. Furthermore,
it is shown there that $\De$ is generated by the matrices of the form
$\left(\matrix A & 0 \\ 0 & \sideset{^t}{^{-1}}{A}\endmatrix\right)$
where $A\in\SL_g(\Z)$,
$\left(\matrix 1 & B \\ 0 & 1\endmatrix\right)$
and
$\left(\matrix 1 & 0 \\ B & 1\endmatrix\right)$
where $B$ is integral symmetric $g\times g$ matrix with even diagonal
(here we use the
standard symplectic basis $e_1,\ldots,e_g,f_1,\ldots f_g$ such that
$(e_i,f_j)=\de_{i,j}$).
We claim that $\De\sub [\Ga_{1,2},\Ga_{1,2}]$. For the proof let us introduce
the relevant elementary matrices following the notation of \cite{Cl} 5.3.1.
Let $\S_{2g}$ be the set of pairs $(i,j)$ where $1\le i,j\le 2g$ which
are not of the form $(2k-1,2k)$ or $(2k,2k-1)$. Then for for every
$(i,j)\in \S_{2g}$ we define an elementary
matrix $E_{ij}$ as follows:
$$E_{2k,2l}=
\left(\matrix 1 & 0 \\ \ga_{k,l} & 1\endmatrix\right),$$
$$E_{2k-1,2l-1}=
\left(\matrix 1 & -\ga_{k,l} \\ 0 & 1\endmatrix\right),$$
$$E_{2k-1,2l}=
\left(\matrix e_{kl} & 0 \\ 0 & e_{lk}^{-1}\endmatrix\right),$$
$$E_{2l,2k-1}=E_{2k-1,2l}$$
where $\ga_{kl}$ has zero $(\a,\b)$-entry unless $(\a,\b)=(k,l)$
or $(\a,\b)=(l,k)$, in the latter case $(\a,\b)$-entry is 1;
$e_{kl}$ for $k\neq l$ is the usual elementary matrix with units on
the diagonal and at $(k,l)$-entry and zeros elsewhere.
Using these matrices one can say that $\De$ is generated by
$E_{2k-1,l}$ with $k\neq l$, $E_{2k,2l}$ and $E_{2k-1,2l-1}$
with $k\neq l$, and $E_{i,i}^2$ for all $1\le i\le 2g$.
It remains to notice that all the matrices $E_{ij}$ with $i\neq j$
belong to $\Ga_{1,2}$ and use the following relations (cf. \cite{Cl} 9.2.13):
\begin{enumerate}
\item $[E_{ij},E_{kl}]=E_{il}$, if $(j,k)\not\in\S_{2g}$, $j$
is even, and $i$, $j$, $k$, and $l$ are distinct,
\item $[E_{ij},E_{ki}]=E_{ii}^2$, if $(j,k)\not\in\S_{2g}$, $j$
is even, and $i$, $j$, and $k$ are distinct.
\end{enumerate}
\ed

\subsection{Case of principally polarized abelian surfaces}
Let $A/S$ be a relative abelian surface, $L$ be a symmetric line
bundle trivialized along the zero section. Assume also that $d=1$
that is $L$ gives a principal polarization. Then
$L\simeq \O(\Th)\ot\pi^*(\pi_*L)$ where $\Th\sub A$ is theta-divisor.
         
\begin{prop} Assume that $S$ is smooth.
Then one has the following equalities in $\Pic(S)$:
$$\det\pi_*\O_{\Th}=\ov{\om}_A,$$
$$5\cdot\ov{\om}_A=\de+\De'(L^2),$$
where $\De'(L^2)=\det\pi_*(L^2)+2\cdot\ov{\om}_A$,
$\de$ is the class of the divisor consisting of points $s\in S$
such that $\Th_s$ is singular,
\end{prop}
          
\Pf . First of all, we note that $\det\pi_*\O_{\Th}=\det\pi_*\om_{\Th}$
by Serre's duality. Now by adjunction we have
$$\om_{\Th}=\O_{\Th}(\Th)\ot\pi^*\ov{\om}_A,$$
which implies the first equality due to triviality of
$\det\pi_*(\O_A(\Th))$ and $\det\pi_*(\O_A)$.
We also deduce that
$$\det\pi_*(\om_{\Th}^2)=\det\pi_*\O_{\Th}(2\Th)+6\cdot\ov{\om}_A.$$
Next, since $L^2\simeq\O(2\Th)\ot\pi^*(\pi_*L)^2$ we obtain that
$$\det\pi_*\O_A(2\Th)=\det\pi_*L^2+4\cdot\ov{\om}_A.$$
The exact sequence
$$0\ra\O_A(\Th)\ra\O_A(2\Th)\ra\O_{\Th}(2\Th)\ra 0$$
shows that $\det\pi_*\O_{\Th}(2\Th)=\det\pi_*\O_A(2\Th)$.
Combining it with the above equalities we get
$$\det\pi_*(\om_{\Th}^2)=\det\pi_*L^2+10\cdot\ov{\om}_A=\De'(L^2)+
8\cdot\ov{\om}_A.$$
On the other hand, since $\Th$ is a stable curve over $S$ we have according
to Mumford's Theorem 5.10 in \cite{Mumst}
$$\det\pi_*(\om_{\Th}^2)=13\cdot\det\pi_*\om_{\Th}-\de=
13\cdot\ov{\om}_A-\de.$$
Comparing this with the previous expression for $\det\pi_*(\om_{\Th}^2)$
we obtain the result.
\ed
         
Let $\ov{\MM}_2$ be the moduli stack of stable curves of genus 2,
$\MM_2$ be the open substack corresponding to smooth curves,
$\MM'_2$ be the substack of $\ov{\MM}_2$ corresponding to
curves which are either smooth or reducible.
The Picard groups of these stacks can be described
as follows (see \cite{Mute}, \cite{Mumst}, \cite{Muen}).
$\Pic(\ov{\MM}_2)$ is isomorphic to $\Z^2$ and is generated by the classes
$\de_0$, $\de_1$ and $\la$, where $\de_0$ (resp. $\de_1$)
is the class of the divisor of
singular irreducible curves (resp. reducible curves),
$\la=\det\pi_*\om_{\CC}$ where $\pi:\CC\ra\ov{\MM}_2$ is the universal curve,
with the only relation
\begin{equation}\label{Mure}
10\cdot\la=\de_0+2\cdot\de_1.
\end{equation}
It follows that $\Pic(\MM'_2)$ is generated by $\la$ and $\de_1$
with the relation $10\cdot\la=2\cdot\de_1$, and $\Pic(\MM_2)$
is generated by $\la$ with the relation $10\cdot\la=0$.
Note that the theta divisors $\Th_s$ are either smooth or reducible,
so in the above situation we get a morphism $f:S\ra\MM'_2$ and
our computation shows that $\De'(L)=f^*(5\cdot\la-\de_1)$.
         
\begin{cor} In the above situation $\De(L^{2n})=0$ for any $n$.
\end{cor}
         
\Pf . For $n=1$ this follows from the triviality of
$10\cdot\la-2\cdot\de_1$   in   $\Pic(\MM'_2)$.  Now  the
triviality of $\De(L^{2n})$ in general follows from Proposition \ref{Ln}.
\ed
      
\begin{rem}
Note that $L^2\simeq L_{\phi}:=(\id,\phi)^*\PP$,
where $\phi:A\ra\hat{A}$ is the polarization corresponding to $L$. Hence,
$\De'(L^2)$ is the pull-back of the line bundle $\De'(L_{\phi})$
over the moduli stack $\AA_2$.
The explicit trivialization  of
$2\cdot\De'(L_{\phi})=10\cdot\ov{\om}-2\cdot\de$   in    the
analytic situation can be found by considering
the  following modular  form  of  weight 5 on Siegel half-space
$\SH_2$ (cf. \cite{Muen})
$$f(Z)=\prod_{a,b\ even}\th
\left[\matrix a\\ b\endmatrix\right](0,Z)$$
where the product is taken over all 10 even theta-characteristics.
Then $f$ defines a section of $\ov{\om}^5$ vanishing
precisely on the locus $\De\sub\SH_2$ corresponding to products of elliptic
curves. It is known that $f$ is a modular form for the
group $\Sp_{4}(\Z)$ with a non-trivial character
$\chi_0:\Sp_4(\Z)\ra\{\pm 1\}$ (such a character is unique,
and is obtained from the sign character on
$\Sp_4(\Z/2\Z)\simeq  S_6$).  Thus,  $f^2$   gives   the
$\Sp_4(\Z)$-equivariant trivialization          of
$\ov{\om}^{10}(-2\De)$  which  descends   to   a
trivialization over $\AA_2$.
This  argument also shows that the element
$\De'(L_{\phi})=5\cdot\ov{\om}-\de\in\Pic(\AA_2)$
is non-trivial. In fact, it generates the torsion subgroup of
$\Pic(\AA_2)$ (cf. \cite{GH}).
Furthermore, one can show that pull-backs
of  $\De'(L_{\phi})$  to  either of  two   irreducible
components of $\wt{\AA}_2$ are non-trivial.
Indeed, it is sufficient to check that the subgroup in
$\Sp_4(\Z/2\Z)\simeq S_6$ preserving a quadratic form
$q$  on  $(\Z/2\Z)^4$,  such  that  $q(x+y)+q(x)+q(y)$ is the
symplectic form, contains an odd permutation.
Recall that the identification of $\Sp_4(\Z/2\Z)$ with $S_6$ is
obtained by considering the action on the set of 6 odd quadratic
forms $q$ as above. Using this it is easy to compute that
the matrix $E_{14}$ (see the proof of Theorem \ref{torsion}), preserving
the standard even form $q_0=x_1y_1+x_2y_2$, corresponds to the
product of three transpositions.
Similarly,  the  matrix  $E_{11}$  preserves  the  odd   form
$x_1y_1+x_2y_2+x_2^2+y_2^2$    and    corresponds    to    a
transposition.
\end{rem}

\vspace{0.5cm}
Deparment of Mathematics, Harvard University, Cambridge, MA 02138
   
{\it E-mail address}: apolish@@math.harvard.edu
                 
\end{document}